\pdfoutput=1
\RequirePackage{ifpdf}
\ifpdf 
\documentclass[pdftex]{sigma}
\else
\documentclass{sigma}
\fi

\numberwithin{equation}{section}
\newtheorem{Theorem}{Theorem}[section]
\newtheorem{Lemma}[Theorem]{Lemma}
\newtheorem{Proposition}[Theorem]{Proposition}
{ \theoremstyle{definition}
\newtheorem{Definition}[Theorem]{Definition}
\newtheorem{Example}[Theorem]{Example}
}

\begin{document}

\def\Hom{\mathop{\rm Hom}\nolimits} \def\Id{\mathop{\rm Id}\nolimits}

\def\cl{\blacktriangleright\hspace{-4pt} <}

\def\al{>\hspace{-4pt}\vartriangleleft}

\def\acl{\blacktriangleright\hspace{-4pt}\vartriangleleft}

\def\build#1_#2^#3{\mathrel{\mathop{\kern 0pt#1}\limits_{#2}^{#3}}}

\newcommand{\ps}[1]{^{(#1)}}

\newcommand{\ns}[1]{_{{\langle #1\rangle}}}

\newcommand{\arXivNumber}{1408.5540}

\allowdisplaybreaks

\renewcommand{\thefootnote}{$\star$}

\renewcommand{\PaperNumber}{093}

\FirstPageHeading

\ShortArticleName{Generalized Coef\/f\/icients for Hopf Cyclic Cohomology}

\ArticleName{Generalized Coef\/f\/icients for Hopf Cyclic Cohomology\footnote{This paper is a~contribution to the
Special Issue on Noncommutative Geometry and Quantum Groups in honor of Marc A.~Rief\/fel.
The full collection is available at
\href{http://www.emis.de/journals/SIGMA/Rieffel.html}{http://www.emis.de/journals/SIGMA/Rieffel.html}}}

\Author{Mohammad HASSANZADEH, Dan KUCEROVSKY and Bahram RANGIPOUR}

\AuthorNameForHeading{M.~Hassanzadeh, D.~Kucerovsky and B.~Rangipour}

\Address{University of New Brunswick, Fredericton, Canada}
\Email{\href{mailto:mhassan@uwindsor.ca}{mhassan@uwindsor.ca},
\href{mailto:dkucerov@unb.ca}{dkucerov@unb.ca},
\href{mailto: bahram@unb.ca}{bahram@unb.ca}}

\ArticleDates{Received August 19, 2013, in f\/inal form August 22, 2014; Published online September 01, 2014}

\Abstract{A category of coef\/f\/icients for Hopf cyclic cohomology is def\/ined.
It is shown that this category has two proper subcategories of which the smallest one is the known category of stable
anti Yetter--Drinfeld modules.
The middle subcategory is comprised of those coef\/f\/icients which satisfy a~generalized SAYD condition depending on both
the Hopf algebra and the (co)algebra in question.
Some examples are introduced to show that these three categories are dif\/ferent.
It is shown that all components of Hopf cyclic cohomology work well with the new coef\/f\/icients we have def\/ined.}

\Classification{19D55; 16T05; 11M55}

\rightline{\it Dedicated to Professor Marc A.~Rieffel}

\renewcommand{\thefootnote}{\arabic{footnote}}
\setcounter{footnote}{0}

\section{Introduction}

A suitable class of coef\/f\/icients for Hopf cyclic cohomology was f\/irst introduced in~\cite{HaKhRaSo1,JS}.
These coef\/f\/icients generalize modular pairs in involution discovered by Connes--Moscovici for Hopf cyclic
cohomology~\cite{ConMos:HopfCyc}, and are called stable anti Yetter--Drinfeld (SAYD) modules.
Hopf cyclic cohomology is def\/ined for a~datum consisting of a~symmetry, a~space, and a~space of coef\/f\/icients.
The symmetry is def\/ined by a~Hopf algebra, the space is def\/ined by an algebra or a~coalgebra on which the Hopf algebra
acts or coacts, and the coef\/f\/icients space is a~module and comodule over the Hopf algebra satisfying two conditions
def\/ined merely by the Hopf algebra structure.
One notes that bialgebra cyclic cohomology introduced in~\cite{Kay} extends the Hopf cyclic cohomo\-lo\-gy tremendously~by
relaxing the symmetry to be a~bialgebra and the coef\/f\/icients to be just stable module-comodule.
However for computing the latter cohomology one cannot rely on the homological algebra of coalgebras and algebras.
What we introduce in this paper is a~class of coef\/f\/icients larger than SAYD modules but still computable in terms of
derived category of (co)algebras.
We show that the bialgebra cyclic cohomology for the new coef\/f\/icients coincide with their Hopf cyclic cohomology.

Based on examples, we observe that Hopf cyclic cohomology works well with some coef\/f\/icients which are not SAYD modules
over Hopf algebras.
We study these modules and prove that there are at least three noticeable categories of such coef\/f\/icients.
The idea is to show that the coef\/f\/icients not only depend on the Hopf algebra in question but that they may also be
related to the (co)algebra on which the Hopf algebra (co)acts.
This observation unties our hands in situations where the Hopf algebra lacks a~large class of SAYD modules.
Existence of such a~Hopf algebra is already shown in~\cite{rs1}, where it is proved that any f\/inite-dimensional SAYD
module over Connes--Moscovici Hopf algebra is a~direct sum of the original modular pair in involution found
in~\cite{ConMos:HopfCyc}.

We cover two cases of Hopf cyclic cohomology in this paper.
For an~$H$-module coalgebra~$C$ we def\/ine $_CH$-SAYD modules and for an~$H$-module algebra~$A$ we introduce $_AH$-SAYD
modules.
We show that SAYD modules are both $_CH$-SAYD modules and $_AH$-SAYD modules for any~$A$ and~$C$, however we give
several examples showing that these inclusions of categories are proper.

Next we show that there are examples of coef\/f\/icients for~$H$-module (co)algebras $(C)A$ which are not even $(_CH)_AH
$-SAYD modules, while the correspondent Hopf cyclic complexes are well-def\/ined.
We call these coef\/f\/icient Hopf cyclic coef\/f\/icients (HCC).
It is not dif\/f\/icult to see that the cup product in Hopf cyclic cohomology works well with these new coef\/f\/icients.

One notes that after introducing Hopf cyclic cohomology with coef\/f\/icients in~\cite{HaKhRaSo1,JS}, there are some
generalizations of the theory in~\cite{bs2,Kay,kk}.
However in all computable examples the authors had to transform their examples into the original Hopf cyclic complexes
of some new Hopf algebras.
As an important example of such, in~\cite{Kay} the coef\/f\/icients are module/comodules over bialgebras and do not have to
satisfy any YD like condition.
However, the author in~\cite{Kay} proves that if these coef\/f\/icients are SAYD, when the bialgebra is Hopf algebra, then
the Hopf cyclic cohomology and bialgebra cyclic cohomoloy coincide.
We prove that the same result holds for new coef\/f\/icients.
The proof is a~modif\/ication of the proof given in~\cite{Kay} for SAYD modules over Hopf algebras.

Hopf cyclic cohomology generalizes Lie algebra cohomology.
One of the most useful features of Lie algebra cohomology is the induction/restriction functor which simply allows us to
use a~subalebra of the Lie algebra in question to calculate the Lie algebra cohomology of the Lie algebra with
coef\/f\/icients in the induced module in terms of the Lie algebra cohomology of the subalgebra with coef\/f\/icients in the
given module.
This feature is successfully adopted for Hopf cyclic cohomology in~\cite{JS}.
Thanks to one of the referees some of our examples can be calculated via the induction procedure.

\textbf{Notations.} In this paper we denote a~Hopf algebra by~$H$ and its counit by~$\varepsilon$.
We use the Sweedler summation notation $\Delta(h)= h\ps{1}\otimes h\ps{2}$ for the coproduct of a~Hopf algebra.
Furthermore $\blacktriangledown(h)= h\ns{-1}\otimes h\ns{0}$ and $\blacktriangledown(h)= h\ns{0}\otimes h\ns{1}$ are
used for the left and right coactions of a~coalgebra, respectively.

\section{Hopf cyclic cohomology coef\/f\/icients}

In this section we introduce two generalized versions of SAYD modules for module coalgebras and module algebras.
Let a~Hopf algebra act on an algebra or a~coalgebra, this action being compatible with (co)algebra structure.
The main idea here is to show that there are coef\/f\/icients by which the Hopf cyclic cohomology of the datum is
well-def\/ined where these coef\/f\/icients are not only (co)representations of the Hopf algebras but also represent the
algebra or coalgebra in question.
Let us recall from~\cite{HaKhRaSo1} that a~right-left SAYD module~$M$ over a~Hopf algebra~$H$ is a~right module and
a~left comodule over~$H$ such that
\begin{gather*}
 \blacktriangledown(mh)= S\big(h^{(3)}\big) m\ns{-1}h^{(1)}\otimes m\ns{0}h^{(2)}
\qquad
\text{(AYD condition)},
\\
 m\ns{0}m\ns{-1}=m
\qquad
\text{(stability condition)}.
\end{gather*}

Since bicrossed product Hopf algebras will be used repeatedly in this paper, let us recall this notion here.
Let ${\cal U}$ and ${\cal F}$ be Hopf algebras and furthermore let ${\cal U}$ be a~right ${\cal F}$-comodule coalgebra.
One then forms a~cocrossed product coalgebra ${\cal U} \cl {\cal F} $ that has ${\cal
F}\otimes {\cal U}$ as underlying vector space and the following coalgebra structure
\begin{gather*}
 \Delta(f\cl u)= f\ps{1}\cl u\ps{1}\ns{0}\otimes f\ps{2}u\ps{1}\ns{1}\cl u\ps{2},
\qquad
 \epsilon(f\cl u)=\epsilon(f)\epsilon(u).
\end{gather*}
Dually if ${\cal F}$ is a~left ${\cal U}$-module algebra, we can endow the underlying vector space ${\cal F}\otimes
{\cal U}$ with an algebra structure, to be denoted by ${\cal F}\al {\cal U}$, with $1\al 1$ as its unit and the product
given~by
\begin{gather*}
(f\al u)(g\al v)=f
u\ps{1} g\al u\ps{2}v.
\end{gather*}
The Hopf algebras ${\cal U}$ and ${\cal F}$ as above, with the action and coaction as above, form a~matched pair of Hopf
algebras if they satisfy the following compatibility conditions for any $u\in{\cal U}$, and any $f\in {\cal F}$
\begin{gather*}
 \epsilon(u f)=\epsilon(u)\epsilon(f),
\qquad
 \Delta(u f)=u\ps{1}\ns{0} f\ps{1}\otimes u\ps{1}\ns{1}\big(u\ps{2} f\ps{2}\big),
\qquad
 \blacktriangledown(1)=1\otimes 1,
\\
 \blacktriangledown(uv)=u\ps{1}\ns{0} v\ns{0}\otimes u\ps{1}\ns{1}\big(u\ps{2} v\ns{1}\big),
\qquad
 u\ps{2}\ns{0}\otimes\big(u\ps{1}f\big)u\ps{2}\ns{1}=u\ps{1}\ns{0}\otimes u\ps{1}\ns{1}\big(u\ps{2} f\big).
\end{gather*}
One forms a~new Hopf algebra ${\cal F}\acl {\cal U}$, called the bicrossed product of the matched pair $({\cal F}, {\cal
U})$; it has ${\cal F}\cl {\cal U}$ as underlying coalgebra, ${\cal F}\al {\cal U}$ as underlying algebra and the
antipode is def\/ined~by
\begin{gather*}
S(f\acl u)=(1\acl S(u\ns{0}))(S(fu\ns{1})\acl 1),
\qquad
f \in {\cal F}, \quad u \in {\cal U}.
\end{gather*}

Throughout the paper all Hopf algebras are assumed to have invertible antipodes.

\subsection[The $_CH$-SAYD modules for module coalgebras]{The $\boldsymbol{{}_CH}$-SAYD modules for module coalgebras}
\label{S1}

For basics of (co)cyclic modules and their correspondence cyclic (co)homology we refer the reader to~\cite{NCG}
and~\cite{Loday}.
If~$M$ is a~right module and left comodule over a~Hopf algebra~$H$, then for any left~$H$-module coalgebra~$C$ one
def\/ines the following para-cocyclic structure on $C^{*}(C, M):=M\otimes C^{\otimes (n+1)} $
\begin{gather}
\partial_i(m\otimes \widetilde{c})= m\otimes c_0\otimes \dots \otimes c_i\ps{1}\otimes c_i\ps{2}\otimes \dots \otimes c_n,\nonumber
\\
\nonumber
 \partial_n(m\otimes \widetilde{c})= m\ns{0}\otimes c_0\ps{2}\otimes c_1\otimes \dots \otimes c_{n}\otimes
m\ns{-1}c_0\ps{1},
\\
 \sigma_i(m\otimes \widetilde{c})=m\otimes c_0\otimes \dots \otimes \varepsilon(c_i)\otimes \dots \otimes c_n,\nonumber
\\
 \tau_n(m\otimes \widetilde{c})=m\ns{0}\otimes c_1\otimes \dots \otimes c_n\otimes m\ns{-1} c_0,\label{module coalgebra}
\end{gather}
where $\widetilde{c}=c_0\otimes \dots \otimes c_n$.
In addition, when~$M$ is SAYD module over~$H$, it is shown in~\cite{HaKhRaSo2} that the above para-cocyclic structure is
well-def\/ined on $C^n_H(C,M):=M\otimes_H\otimes C^{\otimes n+1}$.
Here~$H$ acts on $C^{\otimes n+1}$ diagonally.

However the cyclic structure on $C^{*}_{H}(C, M)$ can be well-def\/ined for some other types of mo\-du\-les~$M$ which are not
SAYD over~$H$.

\begin{Definition}
Let~$C$ be a~left~$H$-module coalgebra.
A~right-left module-comodule~$M$ over~$H$ is called a~$_CH$-SAYD module if for any $m\in M$, $h\in H$, $c\in C$,
$\tilde{d}\in C^{\otimes \bullet}$ we have
\begin{gather*}
\text{\rm ($_CH$-AYD)}
\qquad
(mh)\ns{-1} c \otimes (mh)\ns{0}= S\big(h\ps{3}\big)m\ns{-1}h\ps{1} c \otimes m\ns{0}h\ps{2},
\\
\text{\rm ($_CH$-stability)}
\qquad  m\ns{0}
\otimes_{H} m\ns{-1} \widetilde{d}=m\otimes_{H} \widetilde{d}.
\end{gather*}
Here the action of~$H$ on $C^{\otimes\bullet}$ is diagonal.
We use the symbol $_CH$-$\cal{SAYD}$ for the category whose objects are $_CH$-SAYD modules and whose morphisms
are~$H$-module and~$H$-comodule morphisms.

\end{Definition}

Let us introduce some examples of $_CH$-SAYD modules which are not necessarily a~SAYD module over Hopf algebras.
First we generalize the notion of modular pair in involution~\cite{ConMos:HopfCyc}.
\begin{Lemma}
\label{first}
Let~$C$ be an~$H$-module coalgebra,~$\delta$ a~character and~$\sigma$ a~group-like element for~$H$.
Let $(\delta, \sigma)$ be a~modular pair, i.e.\ $\delta(\sigma)=1$ and $_CH$-in involution, i.e.\
\begin{gather*}
S^2_{\delta}(h) c= \big(\sigma h \sigma^{-1}\big) c,
\qquad
h\in H, \quad c\in C,
\end{gather*}
where $S_{\delta}(h)= \delta\big(h\ps{1}\big)S\big(h\ps{2}\big)$.
Then $^\sigma\mathbb{C}_{\delta}$ is a~$_CH$-{\rm SAYD} module via the action and coaction induced by~$\delta$ and~$\sigma$
respectively.
\end{Lemma}
\begin{proof}
One checks that $S_{\delta}^{-1}(h)= S^{-1}\big(h\ps{1}\big)\delta\big(h\ps{2}\big)$. On the other hand we have
\begin{gather*}
 S\big(h\ps{3}\big)\sigma h\ps{1} c\otimes \delta\big(h\ps{2}\big)= S_{\delta}\big(h\ps{2}\big)\sigma h\ps{1} c \otimes 1
= S^2_{\delta}\big(S_{\delta}^{-1}\big(h\ps{2}\big)\big)\sigma h\ps{1}c\otimes 1\\
\qquad {}= \sigma S_{\delta}^{-1}\big(h\ps{2}\big)\sigma^{-1}\sigma h\ps{1} c\otimes 1
= \sigma S^{-1}\big(h\ps{2}\big)\delta\big(h\ps{3}\big) h\ps{1}c\otimes 1= \sigma c \otimes \delta(h).
\end{gather*}
Using the fact that $(\delta, \sigma)$ is a~modular pair, the $_CH$-stability condition is obvious.
\end{proof}
For any~$H$-module coalgebra~$C$, the following lemma introduce a~one-dimensional $_DH$-SAYD module for some~$H$-module
coalgebra~$D$.
\begin{Lemma}
\label{anysym}
Let $(\delta,\sigma)$ be a~modular pair for the Hopf algebra~$H$ and $C$ be a~left~$H$-module coalgebra.
We define the following subspace of~$C$
\begin{gather*}
I=\big\{S_{\delta}^2(h) c- \big(\sigma h \sigma^{-1}\big) c
\,
\big|
\,
\text{for all}
\;
h\in H,
\;
\text{and}
\;
c\in C\big\}.
\end{gather*}
Then~$I$ is a~coideal of~$C$ and $D:=\frac{C}{I}$ is an~$H$-module coalgebra.
Furthermore $^\sigma\mathbb{C}_{\delta}$ is a~$_DH$-{\rm SAYD} module.
\end{Lemma}
\begin{proof}
Since
\begin{gather*}
S^2_\delta(h)= \delta\big(h\ps{1}\big)S^2\big(h\ps{2}\big)\delta\big(S\big(h\ps{3}\big)\big),
\end{gather*}
is a~coalgebra map, for any $c\in C$ and $h\in H$ we have
\begin{gather*}
 \Delta\left(S_{\delta}^2(h) c- \sigma h \sigma^{-1} c\right)
 = S_{\delta}^2\big(h\ps{1}\big) c\ps{1}\otimes S_{\delta}^2\big(h\ps{2}\big) c\ps{2} - \big(\sigma h\ps{1} \sigma^{-1} c\ps{1}\big)\otimes
\big(\sigma h\ps{2} \sigma^{-1} c\ps{2}\big)
\\
\qquad
 = S_{\delta}^2\big(h\ps{1}\big) c\ps{1}\otimes S_{\delta}^2\big(h\ps{2}\big) c\ps{2} - \sigma h\ps{1} \sigma^{-1} c\ps{1}\otimes
S_{\delta}^2\big(h\ps{2}\big) c\ps{2}
\\
\qquad\phantom{=}
 - \big(\sigma h\ps{1} \sigma^{-1} c\ps{1}\big)\otimes \big(\sigma h\ps{2} \sigma^{-1} c\ps{2}\big)+ \sigma h\ps{1} \sigma^{-1}
c\ps{1}\otimes S_{\delta}^2\big(h\ps{2}\big) c\ps{2}
\\
\qquad
 = \big(S_{\delta}^2\big(h\ps{1}\big) c\ps{1}- \sigma h\ps{1} \sigma^{-1} c\ps{1}\big)\otimes S_{\delta}^2\big(h\ps{2}\big) c\ps{2}
\\
\qquad\phantom{=}
 + \big(\sigma h\ps{1} \sigma^{-1} c\ps{1}\big)\otimes \big(S_{\delta}^2\big(h\ps{2}\big) c\ps{2}- \big(\sigma h\ps{2} \sigma^{-1}c\ps{2}\big)\big).
\end{gather*}
This shows $\Delta(I)\subseteq I\otimes C + C\otimes I$.
Also we have $\varepsilon(y)=0$ for all $y\in I$.
Therefore~$I$ is a~coideal of~$C$ and $D=\frac{C}{I}$ is coalgebra.
Now it is enough to show that~$I$ is an~$H$-module.
To do this we show that for all $h,g \in H$ and $c\in C$ we have
\begin{gather}
\label{kin10}
hS_{\delta}^2(g)c- h\sigma g \sigma^{-1} \in I.
\end{gather}
First we observe that
\begin{gather*}
\sigma (\sigma^{-1}h \sigma) \sigma^{-1}\big(S_{\delta}^2(g)c\big)- S_{\delta}^2\big(\sigma^{-1}h \sigma\big)S_{\delta}^2(g)c\in I.
\end{gather*}
Since $S_{\delta}^2$ is an algebra map, we have
\begin{gather}
\label{kin11}
hS_{\delta}^2(g)c-\sigma^{-1}S_{\delta}^2(h)\sigma S_{\delta}^2(g)c \in I.
\end{gather}
On the other hand we have
\begin{gather*}
S_{\delta}^2\big(\sigma^{-1}h\sigma g\big)c- \sigma\big(\sigma^{-1}h\sigma g\big) \sigma^{-1} \in I.
\end{gather*}
Therefore
\begin{gather}
\label{kin12}
\sigma^{-1}S_{\delta}^2(h)\sigma S_{\delta}^2(g)c- h\sigma g \sigma^{-1} \in I.
\end{gather}
Adding relations~\eqref{kin11} and~\eqref{kin12}, we obtain~\eqref{kin10}.
Therefore~$I$ is an~$H$-module coalgebra and using def\/inition of~$I$, $(\delta, \sigma)$ is a~$_DH$-modular pair in
involution.
Thus by Lemma~\ref{first} we obtain that $^{\sigma}\mathbb{C}_{\delta}$ is a~$_DH$-SAYD module.
\end{proof}

To exemplify, let us recall that ${\mathcal{H}}:=\mathcal{H}_1^{\rm cop}$, the co-opposite Hopf algebra of the
Connes--Moscovici Hopf algebra in codimension one,~\cite{ConMos:HopfCyc} is generated by the elements $X$, $Y$
and $\delta_k$, where $k\in \mathbb{N}$, subject to the following relations
\begin{gather*}
[Y,X]=X,
\qquad
[Y, \delta_k]=k\delta_k,
\qquad
[X,\delta_k]=\delta_{k+1},
\qquad
[\delta_k, \delta_i]=0,
\qquad
k,i\in \mathbb{N}.
\end{gather*}
Its coalgebra structure and antipode are def\/ined as follows
\begin{gather*}
 \Delta(Y)=1\otimes Y+ Y\otimes 1, \Delta(\delta_1)=\delta_1\otimes 1+ 1\otimes \delta_1,
\\
 \Delta(X)=X\otimes 1+1\otimes X+ Y\otimes \delta_1,
\\
 \varepsilon(X)=\varepsilon(Y)=\varepsilon(\delta_k)=0,
\qquad
 S(X)=-X+ Y\delta_1,
\qquad
S(Y)=-Y,
\qquad
S(\delta_k)=-\delta_k.
\end{gather*}
Let $\mathcal{U}$ \looseness=-1 be the enveloping algebra of the Lie algebra generated by~$X$ and~$Y$.
Also suppose $\mathcal{F}$ is the co-opposite Hopf algebra of the Hopf algebra generated by all the $\delta_k$.
It is shown in~\cite{rm} that~$H$ is canonically isomorphic to the bicrossed product Hopf algebra $\mathcal{F}\acl
\mathcal{U}$ where $\mathcal{U}$ acts on $\mathcal{F}$ via
\begin{gather*}
X \delta_k= \delta_{k+1},
\qquad
Y \delta_k=k\delta_k,
\end{gather*}
and $\mathcal{F}$ coacts on $\mathcal{U}$ via
\begin{gather*}
\blacktriangledown(X)= X\otimes 1+ Y\otimes \delta_1,
\qquad
\blacktriangledown(Y)=Y\otimes 1.
\end{gather*}
A~similar decomposition exists for the higher dimensions of this Hopf algebra denoted ${\mathcal{H}}_n^{\rm cop}$.
Its is shown in~\cite[Theorem 3.16]{rs1} that there is only one-dimensional SAYD module for the co-opposite Hopf algebra
of the Connes--Moscovici Hopf algebra ${\mathcal{H}}^{\rm cop}_{1}$, namely $^1\mathbb{C}_{\delta}$.
Since for ${\mathcal{H}}^{\rm cop}_{1}$ we have $S^2\neq \Id$, it follows that $\delta\neq \varepsilon$.
But in the following example we show that $^1\mathbb{C}_{\varepsilon}$ def\/ines a~$_C{\mathcal{H}}^{\rm cop}_{1}$-SAYD
module.

\begin{Example}
In Lemma~\ref{anysym}, let ${\mathcal{H}}= {\mathcal{H}}^{\rm cop}_{1}$, $\delta=\varepsilon$, $\sigma=1$ and $C=
{\mathcal{H}}$, where the ${\mathcal{H}}$-module structure is given by the multiplication.
The relations
\begin{gather*}
S^2(Y)=Y,
\qquad
S^2(\delta_k)=\delta_k,
\qquad
S^2(X)=X-\delta_1,
\end{gather*}
show that~$I$ is the coideal generated by the elements of the form
\begin{gather*}
\big\{S^2(X^n)h- X^n h \,\big| \, n\in \mathbb{N}\big\}.
\end{gather*}
Using $S^2(X)=X-\delta_1$, $X\delta_i= \delta_{i+1}+\delta_i X$ and the fact that $S^2$ is a~morphism of algebras, one
obtains the following relation easily by induction
\begin{gather*}
S^2(X^n)-X^n= \delta_n + \sum\limits_{{i=1}}^{n-1}\delta_i h_{n,i},
\end{gather*}
where $h_{n,i}$ are certain elements in ${\mathcal{H}}$.
This shows that the coideal~$I$ is generated~by
\begin{gather*}
\{\delta_i h,
\;
h\in {\mathcal{H}},
\;
i\in\mathbb{N} \}.
\end{gather*}
We def\/ine a~coalgebra map $\varphi: {\mathcal{H}}\longrightarrow {\cal U}$, given~by
\begin{gather*}
f\acl u\longmapsto \varepsilon(f)u.
\end{gather*}
One checks that $ker(\varphi)=I$ and therefore $\frac{{\mathcal{H}}}{I}\cong {\cal U}$.
As a~result the action of ${\mathcal{H}}$ on $D\cong {\cal U}$ which has been def\/ined by multiplication simplif\/ies to
\begin{gather*}
(f\acl u)(1\acl v)= \varphi(f \acl uv)= \varepsilon(f)uv,
\end{gather*}
where $f\acl u \in {\mathcal{H}}$ and $1\acl v\in {\cal U}\cong D$.
Therefore $D=\frac{{\mathcal{H}}}{I}\cong {\cal U}$ is a~${\mathcal{H}}$-module coalgebra and~$^1\mathbb{C}_{\varepsilon}$ is a~${_{\cal U}}{\mathcal{H}}$-SAYD module.
\end{Example}

\subsection{HCC modules for coalgebras}

\begin{Definition}
Let~$C$ be a~coalgebra which is also an~$H$-module but not necessarily module-coalgebra where the action of~$H$ on
$C^{\otimes \bullet}$ is diagonal.
A~module-comodule~$M$ over~$H$ is called a~$(H, C)$-Hopf cyclic coef\/f\/icients, abbreviated by $_CH$-HCC, if the
cosimplicial and cyclic ope\-ra\-tors on $C^{*}_{H}(C, M)$ are well-def\/ined and turn it to a~cocyclic module.
We use the symbol $_CH$-$\cal{HCC}$ for the category whose objects are $_CH$-HCC modules and whose morphisms
are~$H$-module and~$H$-comodule morphisms
\end{Definition}
In the following statement we show that for any $_CH$-SAYD module the cocyclic module structure def\/ined in~\eqref{module
coalgebra} is well-def\/ined and therefore we obtain a~Hopf cyclic cohomology with these new coef\/f\/icients:
\begin{Proposition}
\label{prop-01}
Let~$C$ be an~$H$-module coalgebra.
Then any $_CH$-{\rm SAYD} module~$M$ is a~$_CH$-{\rm HCC}.
\end{Proposition}
\begin{proof}
Since~$C$ is assumed to be an~$H$-module coalgebra, all cosimplicial structure operators except the very last coface are
well-def\/ined.
On the other hand, the last coface is a~composition of the f\/irst coface and the cyclic operator.
Therefore it is enough to show that the cyclic operator is well-def\/ined.
Let $\pi:C^n(C,M)\rightarrow C^n_H(C,M)$ be the natural projection then we see that the same proof as
in~\cite{HaKhRaSo2} works here
\begin{gather*}
\pi(\tau_n(mh\otimes_{H} c_0\otimes\dots \otimes c_n))
 =(mh)\ns{0}\otimes_H c_1\otimes\dots \otimes c_n\otimes (mh)\ns{-1}c_0
\\
\qquad
 =m\ns{0}h\ps{2}\otimes_{H} c_1\otimes\dots \otimes c_n \otimes S\big(h\ps{3}\big)m\ns{-1}h\ps{1}c_0
\\
\qquad
 =m\ns{0}\otimes_{H} h\ps{2}c_1\otimes \dots \otimes h\ps{n+1}c_n \otimes h\ps{n+2} S\big(h\ps{n+3}\big)m\ns{-1}h\ps{1}c_0
\\
\qquad
 =m\ns{0}\otimes_{H} h\ps{2}c_1\otimes\dots \otimes h\ps{n+1}c_n \otimes m\ns{-1}h\ps{1}c_0
\\
\qquad
 =\pi\big(\tau_n\big(m\otimes h\ps{1}c_0\otimes\dots \otimes h\ps{n+1}c_n\big)\big).
\end{gather*}
We use the $_CH$-AYD module condition in the second equality.
To verify the cocyclicity condition, using the $_CH$-stability condition one has
\begin{gather*}
\tau_n^{n+1}(m\otimes_{H}\widetilde{c})= m\ns{0} \otimes_{H} m\ns{-1}\widetilde{c}=m\otimes_{H}\widetilde{c}.
\tag*{\qed}
\end{gather*}
\renewcommand{\qed}{}
\end{proof}

\begin{Lemma}
Let~$C$ be an~$H$-module coalgebra.
If the action of~$H$ on~$C$ is cocommutative, i.e.\
\begin{gather}
\label{cond1}
h\ps{1}c_1 \otimes h\ps{2} c_2 = h\ps{2} c_1 \otimes h\ps{1} c_2,
\qquad
h\in H,
\quad
c_1, c_2 \in C,
\end{gather}
then any module~$M$ over~$H$, with the trivial coaction, defines a~$_CH$-{\rm HCC}.
\end{Lemma}
\begin{proof}
Since~$C$ is assumed to be an~$H$-module coalgebra it suf\/f\/ices to show that $\tau_n$ is well-def\/ined.
Indeed, let $\pi:C^n(C,M)\rightarrow C^n_H(C,M)$ be the natural projection, then
\begin{gather*}
 \pi\big(\tau_n\big(m\otimes h\ps{1}c_0\otimes \dots \otimes h\ps{n+1} c_n\big)\big)
 =m\otimes_H h\ps{2}c_1\otimes \dots \otimes h\ps{n+1}c_n \otimes h\ps{1} c_0
\\
\qquad{}
 =m\otimes_H h\ps{1}c_1\otimes \dots \otimes h\ps{n}c_n \otimes h\ps{n+1} c_0
 =mh\otimes_H c_1\otimes \dots \otimes c_n\otimes c_0
\\
\qquad{}
 =\pi(\tau_n(mh\otimes_H c_0\otimes \dots \otimes c_n)).
\end{gather*}
We use~\eqref{cond1} in the second equality.
\end{proof}

The following statement is a~dual result to the previous lemma.
\begin{Lemma}
Let~$C$ an~$H$-module coalgebra.
If~$H$ acts on~$C$ commutatively, i.e.\
\begin{gather*}
hgc=ghc,
\qquad
h,g \in H,
\quad
c\in C,
\end{gather*}
then any comodule~$M$ over~$H$, by endowing it with the trivial action, defines a~$_CH$-{\rm SAYD} module.
\end{Lemma}
\begin{proof}
Using the fact that the action of~$H$ on~$M$ is via~$\varepsilon$ and~$H$ acts commutatively on~$C$ we have
\begin{gather*}
S\big(h\ps{3}\big)m\ns{-1} h\ps{1} c\otimes m\ns{0}h\ps{2}
=h\ps{1}S\big(h\ps{2}\big)m\ns{-1}c\otimes m\ns{0}= m\ns{-1}c\otimes m\ns{0}\varepsilon(h).
\end{gather*}
The $_CH$-stability condition is obvious by the triviality of the action.
\end{proof}
\begin{Example}
\label{groupnormal}
Let~$G$ be a~discrete group acting on a~set~$X$ normally, i.e.\ $ghx=hgx$ for all $g,h\in G$ and $x\in X$.
Then consider the group algebra $H=\mathbb{C}G$ and the obvious coalgebra structure $C_X$ on~$X$ where each element is
considered to be a~group-like element.
Obviously $C_X$ is a~module coalgebra over $\mathbb{C}G$ by the action of~$G$ on~$X$.
The action of~$H$ on $C_X$ is commutative.
This introduces a~source of $_CH$-SAYDs which are not ordinary SAYD modules over $H=\mathbb{C}G$.
One notes that all types of SAYD modules~$M$ on $\mathbb{C}G$ are known, see~\cite{HaKhRaSo1}, which in all
cases~$G$-graded vector spaces $M= \bigoplus_g M_g$ where for the (non-trivial) action we have $gm=m$ for all $m\in M_g$
and $hm\in M_{hgh^{-1}}$.
\end{Example}

We present another example.
Since $\mathcal{F}$ is a~Hopf subalgebra of $H={\cal F}\acl {\cal U}$, the Hopf algebra~$H$ acts on the coalgebra
$C:=H\otimes_{\mathcal{F}} \mathbb{C}$ via left multiplication.
It is easy to see that as a~coalgebra~$C$ is isomorphic to $\mathcal{U}$ via $u\longmapsto (1\acl
u)\otimes_{\mathcal{F}}1$.
Via this identif\/ication, the action of~$H$ on~$C$ is as follows.
Let $h=f \acl u \in H$, and $v\in C=\mathcal{U}$.
Then
\begin{gather*}
 (f\acl u)v=((f\acl u)(1\acl v)\otimes_{\mathcal{F}}1)
\\
\hphantom{(f\acl u)v}{}
 =(f\acl uv)\otimes_{\mathcal{F}}1=\big[\big(1\acl u\ps{1}v\ps{1}\big)\big(S\big(\big(u\ps{2}v\ps{2}\big)\big) f\acl 1\big)\big]\otimes_{\mathcal{F}}1
\\
\hphantom{(f\acl u)v}{}
 =\big(1 \acl u\ps{1}v\ps{1}\big)\otimes_{\mathcal{F}}\varepsilon\big(S\big(\big(u\ps{2}v\ps{2}\big)\big) f\big)
 =(1\acl uv)\otimes_{\mathcal{F}}\varepsilon(f)=\varepsilon(f)uv.
\end{gather*}
One easily checks that for any bicrossed product Hopf algebra ${\cal F} \acl {\cal U}$, the Hopf algebra ${\cal U}$ is
an~$H$-module coalgebra via this action.

\begin{Example}
Consider any bicrossed product Hopf algebra $H={\cal F} \acl {\cal U}$ where ${\cal U}$ is a~cocommutative Hopf
algebra.
We have seen that $C={\cal U}$ is an~$H$-module coalgebra by the following action
\begin{gather}
\label{normal-com}
(f\acl u)v= \varepsilon(f)uv.
\end{gather}
Since ${\cal U}$ is cocommutative, it is obvious that ${\mathcal{H}}$ acts on~$C$ cocommutatively.
Therefore any comodule~$M$ over~$H$ with the trivial action def\/ines a~$_{\cal U} H$-SAYD module.
\end{Example}

\begin{Example}
\label{cocom-action}
Consider any bicrossed product Hopf algebra $H={\cal F} \acl {\cal U}$ where ${\cal U}$ is a~cocommutative Hopf
algebra.
Consider the module coalgebra ${\cal U}$ with the action def\/ined in~\eqref{normal-com}.
We show that~$H$ acts on ${\cal U}$ cocommutatively.
First let us recall that
\begin{gather*}
\Delta(f\acl u)=f\ps{1}\acl u\ps{1}\ns{0}\otimes f\ps{2}u\ps{1}\ns{1}\acl u\ns{2}.
\end{gather*}
Using the cocommutativity of $\mathcal{U}$, for $h:=f\acl u$ we get
\begin{gather*}
 h\ps{1}v^1\otimes h\ps{2}v^2=\varepsilon\big(f\ps{1}\big)u\ps{1}\ns{0}v^1\otimes \varepsilon\big(f\ps{2}u\ps{1}\ns{1}\big)u\ps{2}v^2,
\\
 \varepsilon(f)u\ps{1}v^1\otimes u\ps{2}v^2=\varepsilon(f)u\ps{2}v^2\otimes u\ps{1}v^1= h\ps{2}v^1\otimes h\ps{1}v^2.
\end{gather*}
Therefore any module~$M$ over~$H$, providing it with the trivial coaction of~$H$, def\/ines a~$_{\cal U} H$-HCC.
As an example the co-opposite Hopf algebra of Connes--Moscovici Hopf algebra $\mathcal{H}_1^{\rm cop}\cong \mathcal{F}\acl
\mathcal{U}$ acts cocommutatively on ${\cal U}$.
\end{Example}

Here we introduce an example of $_CH$-HHC which is not a~$_CH$-SAYD module.
We have shown that if~$H$ acts on~$C$ cocommutatively then any module~$M$ def\/ines an $_CH$-HHC.
Indeed in this example, we introduce a~triple $(H,C,M)$ such that~$H$ acts on~$C$ cocommutatively but~$M$ is not
a~$_CH$-SAYD module.

\begin{Example}
Consider ${\mathcal{H}}=\mathcal{H}_1^{\rm cop}\cong \mathcal{F}\acl \mathcal{U}$.
We have seen that $C={\cal U}$ is a~${\mathcal{H}}$-module coalgebra by the action $(f\acl u)v= \varepsilon(f)uv$ and
that this action is a~cocommutative action by Example~\ref{cocom-action}.
Now let $M={\mathcal{H}}$ be a~right ${\mathcal{H}}$-module via multiplication and left ${\mathcal{H}}$-comodule~by
trivial coaction.
Since the action is cocommutative, it is a~$_{\cal U}{\mathcal{H}}$-HCC module.
We show that~$M$ is not a~$_C{\mathcal{H}}$-SAYD.
Let $f\acl u:=\delta_1\otimes X\in {\mathcal{H}}$,
$c:=X\in \mathcal{U}$, and $m:= 1\acl 1\in M$.
First we notice that
\begin{gather}
\nonumber
 (f\acl u)\ps{1}\otimes (f\acl u)\ps{2}\otimes (f\acl u)\ps{3}
\\
\qquad
 =f\ps{1}\acl u\ps{1}\ns{0}\otimes f\ps{2}u\ps{1}\ns{1}\acl u\ps{2}\ns{0}\otimes f\ps{3}u\ps{1}\ns{2}u\ps{2}\ns{1}\acl u\ps{3}.
\label{fu}
\end{gather}
After substituting into the AYD condition we observe that
\begin{gather*}
 S\big( (f\acl u )\ps{3}\big)(f\acl u)\ps{1}c\otimes m(f\acl u)\ps{2}
\\
\qquad
 =S\big(f\ps{3}u\ps{1}\ns{2}u\ps{2}\ns{1}\acl u\ps{3}\big)\big(f\ps{1}\acl u\ps{1}\ns{0}\big)c\otimes mf\ps{2}u\ps{1}\ns{1}\acl u\ps{2}\ns{0}
\\
\qquad
 =S\big(f\ps{2}u\ps{1}\ns{2}u\ps{2}\ns{1}\acl u\ps{3}\big) u\ps{1}\ns{0}c\otimes mf\ps{1}u\ps{1}\ns{1}\acl u\ps{2}\ns{0}
\\
\qquad
 =\big(1\acl S\big(u\ps{3}\ns{0}\big)\big)\big(S\big(f\ps{2}u\ps{1}\ns{2}u\ps{2}\ns{1} u\ps{3}\ns{1}\big)\acl 1\big) u\ps{1}\ns{0}c
 \otimes mf\ps{1}u\ps{1}\ns{1}\acl u\ps{2}\ns{0}
\\
\qquad
 =S\big(X\ps{3}\big)X\ps{1}\ns{0}X\otimes (1\acl 1)\big(\delta_1 X\ps{1}\ns{1}\acl X\ps{2}\big)
\\
\qquad
 =X\ns{0}X\otimes (\delta_1 X\ns{1}\acl 1)+ X\otimes (\delta_1 \acl X)+ S(X)X\otimes (\delta_1 \acl 1)
\\
\qquad
 =X^2\otimes (\delta_1 \acl 1)+YX\otimes \big(\delta_1^2 \acl 1\big)+X\otimes (\delta_1 \acl X)-X^2\otimes (\delta_1 \acl 1)
\\
\qquad
 =X\otimes (\delta_1 \acl X)+ YX\otimes \big(\delta_1^2 \acl 1\big)
 \neq X\otimes (\delta_1\acl X)=c\otimes mh.
\end{gather*}
\end{Example}

\subsection{Relation to Kaygun's cyclic cohomology of bialgebras}

In this subsection we f\/irst recall the bialgebra cyclic cohomology developed by Kaygun in~\cite{Kay} and then show that
for any $_CH$-SAYD module the result of these two cohomology coincide.

Let~$H$ be a~bialgebra,~$C$ be an~$H$-module coalgebra,~$M$ be a~right module and left comodule over~$H$.

We let~$H$ act on $C^n(C,M)$ from left~by
\begin{gather*}
L_g(m\otimes c_0\otimes\dots\otimes c_n)= mS\big(g\ps{1}\big)\otimes g\ps{2} c_0\otimes\dots\otimes g\ps{n+2}c_n,
\end{gather*}
and def\/ines the following subspace of $C^n(H,M)$
\begin{gather*}
W^n(C,M):=\big\{\big[L_g,\tau^i\big](C^n(H,M))\,|\, g\in H, i\in{\mathbb Z}^+\big\}.
\end{gather*}
Here $[L_g,\tau]:= L_g\tau-\tau L_g$.

Let us see that~$W$ is stable under cyclic structure of $C^\ast(C,M)$.
Indeed for~$\tau$ we simply see that
\begin{gather*}
\tau \big[L_g,\tau^i\big]=[\tau,L_g] \tau^i+ \big[L_g,\tau^{i+1}\big].
\end{gather*}
For $\partial_m$, where $0\le m\le n$, we use the identity $\partial_m \tau=\tau \partial_{m+1}$ and the fact that
$\partial_m L_g=L_g\partial_m$ to see that
\begin{gather*}
\partial_m\big[L_g,\tau^i\big]=[L_g,\tau]\partial_{m+1}.
\end{gather*}
For $\partial_{n+1}=\tau\partial_0$ we have
\begin{gather*}
\partial_0\big[L_g,\tau^i\big]= \tau\big[L_g,\tau^i\big]\partial_1=[\tau,L_g] \tau^i\partial_1+ \big[L_g,\tau^{i+1}\big]\partial_1.
\end{gather*}

Using the facts that $L_g\sigma_j=\sigma_jL_g$ and $\tau\sigma_j= \sigma_{j-1}\tau$, we see that for $1\le j\le n$, we
have
\begin{gather*}
\sigma_j\big[L_g,\tau^i\big]= \big[L_g,\tau^i\big]\sigma_{j-1}.
\end{gather*}
Finally for $\sigma_0=\sigma_n\tau$ we obtain
\begin{gather*}
\sigma_0\big[L_g,\tau^i\big]= \sigma_n\tau\big[L_g,\tau^i\big]= [\tau,L_g] \sigma_{n-1}\tau^i+ \big[L_g,\tau^{i+1}\big]\sigma_{n-1}.
\end{gather*}

So $\frac {C^\ast(C,M)}{W^\ast(C,M)}$ is a~paracocyclic~$H$-module.
This paracocyclic module is denoted by ${\mathbb P}{\mathbb C}{\mathbb M}^\ast(C,M)$ in~\cite{Kay}.
On the other hand if~$M$ is stable then
\begin{gather*}
\tau^{n+1}(m\otimes c_0\otimes\dots\otimes c_n)= m\ns{0}\otimes m\ns{-n-1}c_0\otimes\dots\otimes m\ns{-1}c_n
\\
\qquad{}
 =m\ns{0} m\ns{-n-3}S(m\ns{-n-2})\otimes m\ns{-n-1}c_0\otimes\dots\otimes m\ns{-1}c_n
\\
\qquad{}
 =L_{m\ns{-2}}(m\ns{0} m\ns{-1}\otimes c_0\otimes\dots\otimes c_n)
 =L_{m\ns{-1}}(m\ns{0}\otimes c_0\otimes\dots\otimes c_n),
\end{gather*}
which proves that the paracocyclic module ${\mathbb C}{\mathbb M}^\ast(C,M):={\mathbb C}\otimes_H{\mathbb P}{\mathbb
C}{\mathbb M}^\ast(C,M)$ is in fact a~cocyclic module.

\begin{Proposition}
\label{prop-02}
Let~$H$ be a~Hopf algebra,~$C$ be an~$H$-module coalgebra, and~$M$ be a~$_CH$-{\rm SAYD} module.
Then $W\subseteq \ker\pi$, where $\pi:C^n(C,M)\rightarrow C^n_H(C,M)$ is the canonical projection.
\end{Proposition}
\begin{proof}
Let us f\/ix $m\otimes c_0\otimes\dots\otimes c_n\in C^n(C,M)$.
Using the $_CH$-{\rm SAYD} property of~$M$, we see that
\begin{gather*}
 \tau L_g(m\otimes c_0\otimes\dots\otimes c_n)= \tau\big(mS\big(g\ps{1}\big)\otimes g\ps{2}c_0\otimes\dots\otimes g\ps{n+2} c_n\big)
\\
\qquad
 =\big(mS\big(g\ps{1}\big)\big)\ns{0}\otimes g\ps{3}c_1\otimes\dots\otimes g\ps{n+2} c_n\otimes\big(mS\big(g\ps{1}\big)\big)\ns{-1} g\ps{2}c_0
\\
\qquad
 = m\ns{0} S\big(g\ps{2}\big)\otimes g\ps{5}c_1\otimes\dots\otimes g\ps{n+2}c_n\otimes S^2\big(g\ps{1}\big)m\ns{-1} S\big(g\ps{3}\big) g\ps{4}c_0
\\
\qquad
 = m\ns{0} S\big(g\ps{2}\big)\otimes g\ps{3}c_1\otimes\dots\otimes g\ps{n+2}c_n\otimes S^2\big(g\ps{1}\big)m\ns{-1}c_0.
\end{gather*}
On the other hand we have{\samepage
\begin{gather*}
 L_g\tau(m\otimes c_0\otimes\dots\otimes c_n)= L_g(m\ns{0}\otimes c_1\otimes\dots\otimes c_n\otimes m\ns{-1}c_0)
\\
\qquad
 = m\ns{0} S\big(g\ps{1}\big)\otimes g\ps{2} c_1\otimes\dots\otimes g\ps{n+1} c_n\otimes g\ps{n+2}m\ns{-1}c_0.
\end{gather*}
So $\pi(\tau L_g(m\otimes c_0\otimes\dots\otimes c_n))= \varepsilon(g) \tau(m\otimes c_0\otimes\dots\otimes c_n)= \pi(L_g\tau(m\otimes
c_0\otimes\dots\otimes c_n))$.}

To f\/inish the proof one uses the facts that, by Proposition~\ref{prop-01},~$\tau$ on $C_H(C,M)$ is well-def\/ined, and
that $[L_g, -]$ is a~derivation.
Indeed,
\begin{gather*}
 \pi\big(\big[L_g,\tau^i\big](m\otimes c_0\otimes\dots\otimes c_n)\big)
\\
\qquad
 =\pi\big([L_g,\tau]\tau^{i-1}(m\otimes c_0\otimes\dots\otimes c_n)\big) +\tau^{i-1}\pi([L_g,\tau](m\otimes c_0\otimes\dots\otimes c_n))=0.
\tag*{\qed}
\end{gather*}
\renewcommand{\qed}{}
\end{proof}
Finally we prove the main result of this subsection.
\begin{Proposition}
Let~$H$ be a~Hopf algebra,~$C$ be an~$H$-module coalgebra, and~$M$ be a~$_CH$-{\rm SAYD} module.
Then ${\mathbb C}{\mathbb M}^\ast(C,M)$ and $C^\ast_H(C,M)$ are isomorphic as cocyclic modules.
\end{Proposition}

\begin{proof}
The proof is similar to the case of of the SAYD modules~\cite{Kay}.
By Proposition~\ref{prop-02} we have the~$H$-linear map between paracocyclic modules
\begin{gather*}
\pi: \ {\mathbb P}{\mathbb C}{\mathbb M}^n(C,M)\longrightarrow C_H^n(H,M).
\end{gather*}
Here we assume that the action of~$H$ on $C^n_H(C,M)$ is trivial.
It is easily seen that
\begin{gather*}
\partial(L_g(m\otimes c_0\otimes\dots\otimes c_n))= \varepsilon(g)\pi(m\otimes c_0\otimes\dots\otimes c_n).
\end{gather*}
So~$\pi$ induces a~map of cocyclic modules
\begin{gather*}
\Pi: \ {\mathbb C}{\mathbb M}^n(C,M)\longrightarrow C_H^n(C,M).
\end{gather*}
On the other hand since $W^\ast$ is a~para-cocyclic submodule of $C^\ast(C,M)$ we have projections of para-cocyclic
modules.
\begin{gather*}
\pi': \ C^n(C,M)\longrightarrow {\mathbb C}{\mathbb M}^n(C,M).
\end{gather*}
Let us check that $\pi'$ is balanced over~$H$.
In the following we use $[x]$ as the class of $x\in C^n(C,M)$ in ${\mathbb P}{\mathbb C}{\mathbb M}^n(C,M)$ and ${\bf
1}\in {\mathbb C}$ as the generator of ${\mathbb C}$ as trivial~$H$-module.
\begin{gather*}
 \pi'(mg\otimes c_0\otimes\dots\otimes c_n)= {\bf 1}\otimes_H [mg\otimes c_0\otimes\dots\otimes c_n]
 ={\bf 1}\otimes_H \big[mg\ps{1}\varepsilon\big(g\ps{2}\big)\otimes c_0\otimes\dots\otimes c_n\big]\\
\qquad{}
 ={\bf 1}\otimes_H \big[mg\ps{1}S\big(h\ps{2}\big)\otimes h\ps{3}c_0\otimes\dots\otimes h\ps{n+3}c_n\big]
\\
\qquad{}
 ={\bf 1}\otimes_H \big[m\otimes g\ps{1} c_0\otimes\dots\otimes g\ps{n+1}c_n\big]
 = \pi'(m\otimes g\ps{1} c_0\otimes\dots\otimes g\ps{n+1}c_n).
\end{gather*}
So $\pi'$ induces a~map of cocyclic modules
\begin{gather*}
\Pi': \ C_H^n(C,M)\longrightarrow {\mathbb C}{\mathbb M}^n(C,M).
\end{gather*}
It is obvious that~$\Pi$ and $\Pi'$ are inverse to each other.
\end{proof}

\subsection[The $_AH$-SAYD and HCC modules for module algebras]{The $\boldsymbol{{}_AH}$-SAYD and HCC modules for module algebras}

Let~$A$ be a~left~$H$-module algebra and~$M$ be a~right-left SAYD module over~$H$.
We notice that $M\otimes A^{\otimes(n+1)}$ is a~right~$H$-module via $(m\otimes \widetilde{a}) h:= mh\ps{1}\otimes
S(h\ps{2})\widetilde{a}$ and the right~$H$-module structure of the ground f\/ield $\mathbb{C}$ is given by $r h:=
\varepsilon(h)r$ for all $r\in \mathbb{C}$ and $h\in H$.
It is shown that the following cocyclic structure on the space of right~$H$-linear homomorphisms
$C^{n}_{H}(A, M):= \Hom_{H}(M\otimes A^{\otimes (n+1)}, \mathbb{C})$ is well-def\/ined~\cite{HaKhRaSo1}
\begin{gather}
 (d_i f)(m\otimes \widetilde{a})= f(m\otimes a_0\otimes \dots \otimes
 a_i a_{i+1} \otimes \dots \otimes a_n),
\qquad
0\leq i< n,\nonumber
\\
 (d_n f)(m\otimes \widetilde{a})=f\big(m\ns{0}\otimes \big(S^{-1}\big(m\ns{-1}\big)a_n\big)a_0\otimes a_1\otimes \dots \otimes a_{n-1}\big),
\nonumber\\
 (s_n f)(m\otimes \widetilde{a})=f(m\otimes a_0 \otimes\dots\otimes a_i\otimes 1\otimes \dots \otimes a_n),
\qquad
0\leq i\leq n,\nonumber
\\
 (t_n f)(m\otimes \widetilde{a})=f(m\ns{0}\otimes S^{-1}(m\ns{-1})a_n\otimes a_0\otimes \dots \otimes a_{n-1}).\label{module algebra}
\end{gather}
where $\widetilde{a}=a_0\otimes \dots \otimes a_n$.

\begin{Definition}
Let~$A$ be a~left~$H$-module algebra.
A~right-left module-comodule~$M$ over~$H$ is called an $_A H$-SAYD module if for all $m\in M$, $h\in H$, $a\in A$ and
$\varphi \in C^{*}_{H}(A, M)$
\begin{enumerate}\itemsep=0pt
\item[i)] $S^{-1}((mh)\ns{-1}) a\otimes (mh)\ns{0}= S^{-1}\big(m\ns{-1}h\ps{1}\big)h\ps{3} a\otimes
m\ns{0}h\ps{2}$,
\item[ii)] $\varphi(m\ns{0} \otimes S^{-1}(m\ns{-1}) a_0\otimes \dots \otimes a_n)=\varphi(m\otimes \widetilde{a})$.
\end{enumerate}
Here the action of~$H$ on $A^{\otimes (n+1)}$ is diagonal.
We use the symbol $_AH$-$\cal{SAYD}$ for the category whose objects are $_AH$-SAYD modules and whose morphisms
are~$H$-module and~$H$-comodule morphisms
\end{Definition}

\begin{Lemma}
Let~$A$ be an~$H$-module algebra.
Then any {\rm SAYD} module over~$H$ is an $_AH$-{\rm SAYD} module.
\end{Lemma}

\begin{proof}
One easily checks that the following properties are satisf\/ied:
i)~$m\ns{0} S^{-1}(m\ns{-1})=m$,
ii)~$\varphi(m\otimes h \widetilde{a})\!=\! \varphi(mh\otimes \widetilde{a})$ if and only if $\varphi \!\in\! C^{*}_{H}(A, M)$,
i.e.\ $\varphi(mh\ps{1}\!\otimes S(h\ps{2}) \widetilde{a})= \varepsilon(h)\varphi(m\otimes \widetilde{a})$.
\end{proof}

\begin{Definition}\looseness=-1
Let~$A$ be an algebra and an~$H$-module (not necessarily an~$H$-module algebra).
A~module-comodule~$M$ over~$H$ is called an $_AH$-Hopf cyclic coef\/f\/icients ($_AH$-HCC), if the cosimplicial and cyclic
operators on $C^{*}_{H}(A, M)$ are well-def\/ined and turn it into a~cocyclic module.
Here the action of~$H$ on $A^{\otimes (n+1)}$ is diagonal.
We use the symbol $_AH$-$\cal{HCC}$ for the category whose objects are $_AH$-HCC modules and whose morphisms
are~$H$-module and~$H$-comodule morphisms.
\end{Definition}

\begin{Proposition}
Let~$A$ be a~left~$H$-module algebra.
Then any $_AH$-{\rm SAYD} module is an $_AH$-{\rm HCC}.
\end{Proposition}
\begin{proof}
Let~$M$ be a~right-left $_AH$-SAYD module.
To show that~$M$ is $_AH$-HCC it is enough to check that the cyclic map is well-def\/ined
\begin{gather*}
 (t_n f)((m\otimes a_0\otimes \dots \otimes a_n) h)
 =(t_n f)\big(m h\ps{1}\otimes S\big(h\ps{2}\big) (a_0 \otimes \dots \otimes a_n)\big)
\\
 =(t_n f)\big(mh\ps{1}\otimes S\big(h\ps{n+2}\big)a_0\otimes \dots \otimes S\big(h\ps{2}\big)a_{n}\big)
\\
 =f\big(\big(mh\ps{1}\big)\ns{0}\otimes S^{-1}\big(\big(mh\ps{1}\big)\ns{-1}\big)S\big(h\ps{2}\big)a_n\otimes S\big(h\ps{n+3}\big)a_0\otimes \dots
\otimes S\big(h\ps{3}\big)a_{n-1} \big)
\\
 =f\big(m\ns{0}h\ps{2}\otimes S^{-1}\big(S\big(h\ps{3}\big)m\ns{-1}h\ps{1}\big)S\big(h\ps{4}\big)a_n\otimes S\big(h\ps{n+4}\big)a_0\otimes \dots
\otimes S\big(h\ps{5}\big)a_{n-1} \big)
\\
 =f\big(m\ns{0}h\ps{2}\otimes S^{-1}\big(h\ps{1}\big)S^{-1}(m\ns{-1})h\ps{3}S\big(h\ps{4}\big)a_n\otimes S\big(h\ps{n+4}\big)a_0\otimes \dots
\otimes S\big(h\ps{5}\big)a_{n-1} \big)
\\
 =f\big(m\ns{0} \otimes h\ps{2} \big[S^{-1}\big(h\ps{1}\big) S^{-1}(m\ns{-1})a_n\otimes S\big(h\ps{n+2}\big)a_0\otimes \dots \otimes
S\big(h\ps{3}\big)a_{n-1} \big] \big)
\\
 =f\big(m\ns{0} \otimes h\ps{2} S^{-1}\big(h\ps{1}\big) S^{-1}(m\ns{-1})a_n\otimes h\ps{3} S\big(h\ps{2n+2}\big)a_0\otimes \dots
\otimes h\ps{n+2}S\big(h\ps{n+3}\big)a_{n-1} \big)
\\
 =f(m\ns{0}\otimes S^{-1}(m\ns{-1})a_n\otimes a_0\otimes \dots \otimes a_{n-1})\varepsilon(h)
\\
 =(t_n f)(m\otimes a_0\otimes \dots \otimes a_n)\varepsilon(h)
 =\left((t_n f)(m\otimes a_0\otimes \dots \otimes a_n)\right) h.
\end{gather*}
To prove the cyclicity, using the $_AH$-stability condition we have
\begin{gather*}
 \big({t_n}^{n+1} f\big)(m\otimes a_0\otimes \dots \otimes a_n)= (t_nf)\big(m\ns{0}\otimes S^{-1}(m\ns{-1})a_1\otimes \dots \otimes
a_n\otimes a_0\big)
\\
\qquad
 =(t_nf)(m\otimes a_1\otimes \dots \otimes a_n\otimes a_0)=f(m\otimes a_0\otimes \dots \otimes a_n).\tag*{\qed}
\end{gather*}
\renewcommand{\qed}{}
\end{proof}

\begin{Lemma}
Let~$A$ be a~left~$H$-module algebra, and $(\delta, \sigma)$ be a~modular pair and $_AH$-in involution, i.e.\
\begin{gather*}
S^{-2}_{\delta}(h) a= \sigma^{-1} h \sigma a,
\qquad
\text{for any}
\quad
a\in A,
\quad
h\in H.
\end{gather*}
Then $^\sigma\mathbb{C}_{\delta}$ is an $_AH$-{\rm SAYD} module.
\end{Lemma}
\begin{proof}
The following computations show the $_AH$-AYD condition
\begin{gather*}
 S^{-1}\big((1_A h)\ns{-1}\big) a\otimes (1_A h)\ns{0}
 =\delta(h)S^{-1}(\sigma) a\otimes 1= \sigma^{-1} \delta\big(h\ps{1}\big)S\big(h\ps{2}\big)h\ps{3} a\otimes 1
\\
\qquad
 =\sigma^{-1} S_{\delta}\big(h\ps{1}\big)\sigma \sigma^{-1} h\ps{2} a\otimes 1
 =S^{-2}_{\delta}\big(S_{\delta}\big(h\ps{1}\big)\big)\sigma^{-1} h\ps{2} a\otimes 1
\\
\qquad
 =S^{-1}\big(h\ps{1}\big) \sigma^{-1} h\ps{3} a\otimes \delta\big(h\ps{2}\big)
 =S^{-1}\big(1\ns{-1}h\ps{1}\big) h\ps{3} a\otimes 1\ns{0} h\ps{2}.
\end{gather*}
The stability condition is obvious from the modular pair condition.
\end{proof}

\begin{Proposition}
Let $(\delta,\sigma)$ be a~modular pair for~$H$ and $A$ be a~left~$H$-module algebra.
We define the following subspace of~$A$
\begin{gather*}
B=\big\{a\in A
\,
\big|
\,
S^{-2}_{\delta}(h) a= \sigma^{-1} h\sigma a,
\;
\text{for all}
\;
h\in H \big\}.
\end{gather*}
Then~$B$ is an~$H$-module subalgebra of~$A$ and $(\delta,\sigma)$ is a~$_BH$-modular pair in involution.
\end{Proposition}
\begin{proof}
It is obvious that $1\in B$.
Since~$A$ is an~$H$-module algebra,~$\sigma$ is a~group-like element and $S^{-2}_{\delta}$ is a~coalgebra map, the
following computation shows that~$B$ is a~subalgebra of~$A$
\begin{gather*}
 S^{-2}_{\delta}(h)(ab)=\big(S^{-2}_{\delta}\big(h\ps{1}\big) a\big)\big(S^{-2}_{\delta}\big(h\ps{2}\big) b\big)
 = \big(\sigma^{-1} h\ps{1}\sigma a\big)\big(\sigma^{-1} h\ps{2}\sigma b\big)=\sigma^{-1} h\sigma (ab).
\end{gather*}

Let $h\in H$ and $b \in B$.
The following computation shows that $hb\in B$ and therefore that~$B$ is a~left~$H$-module
\begin{gather*}
  \sigma^{-1}k \sigma (hb)= \big(\sigma^{-1}(k\sigma h \sigma^{-1}) \sigma\big)b= S^{-2}_{\delta}(k\sigma
h\sigma^{-1})b
\\
\hphantom{\sigma^{-1}k \sigma (hb)}{}
 = S^{-2}_{\delta}(k)\big[S^{-2}_{\delta}(\sigma h\sigma^{-1})b\big]=S^{-2}_{\delta}(k)\sigma^{-1} \big(\sigma h
\sigma^{-1}\big) \sigma b= S^{-2}_{\delta}(k)(hb).
\end{gather*}
Therefore~$B$ is an~$H$-module subalgebra of~$A$ and $(\delta,\sigma)$ is a~$_BH$-modular pair in involution.
\end{proof}

\begin{Lemma}
Let~$A$ be a~left~$H$-module algebra.
If the action of~$H$ on~$A$ is cocommutative, i.e.\
\begin{gather}
\label{cond3}
h\ps{1} a_1 \otimes h\ps{2} a_2 = h\ps{2} a_1 \otimes h\ps{1} a_2,
\qquad
h\in H, \quad a_1, a_2 \in A,
\end{gather}
then any module~$M$ over~$H$, with the trivial coaction of~$H$, defines an $_AH$-{\rm HCC}.
\end{Lemma}
\begin{proof}
The $_AH$-stability condition is obvious by the triviality of the coaction.
Since~$A$ is assumed to be~$H$-module algebra it suf\/f\/ices to show that the cyclic map~$\tau$ is well-def\/ined
\begin{gather*}
 (t_n f) ((m\otimes a_0\otimes \dots \otimes a_n) h )
 =(t_n f)\big(mh\ps{1}\otimes S\big(h\ps{n+1}\big)a_0\otimes \dots \otimes S\big(h\ps{2}\big)a_n\big)
\\
\qquad
 =f\big(mh\ps{1}\otimes S\big(h\ps{2}\big)a_n\otimes S\big(h\ps{n+1}\big)a_0\otimes \dots \otimes S\big(h\ps{3}\big)a_{n-1}\big)
\\
\qquad
 =f\big(mh\ps{1}\otimes S\big(h\ps{n+1}\big)a_n\otimes S\big(h\ps{n}\big)a_0\otimes \dots \otimes S\big(h\ps{2}\big)a_{n-1}\big)
\\
\qquad
 =f(m\otimes a_n\otimes a_0\otimes \dots \otimes a_{n-1})\varepsilon(h)
 =(t_n f)(m\otimes a_0\otimes \dots \otimes a_n)\varepsilon(h).
\end{gather*}
We use~\eqref{cond3} in the third equality.
\end{proof}

\begin{Lemma}
Let~$A$ be a~left~$H$-module algebra.
If~$H$ acts on~$A$ commutatively, i.e.\
\begin{gather*}
hga=gha,
\qquad
h,g \in H, \quad a\in A,
\end{gather*}
then any comodule~$M$ over~$H$, endowed with the trivial action from~$H$, defines an $_AH$-{\rm SAYD} module.
\end{Lemma}
\begin{proof}
The $_AH$-stability condition is obvious.
The following computation proves the $_AH$-AYD condition
\begin{gather*}
 \varepsilon(h)S^{-1} (m\ns{-1} ) a\otimes m\ns{0}
 =h\ps{2}S^{-1}\big(h\ps{1}\big)S^{-1}(m\ns{-1})a\otimes m\ns{0}
\\
\qquad
 =S^{-1}\big(m\ns{-1}h\ps{1}\big)h\ps{2} a\otimes m\ns{0}
 =S^{-1}\big(m\ns{-1}h\ps{1}\big)h\ps{3} a\otimes m\ns{0}h\ps{2}.
\end{gather*}
We use the commutativity of the action on the second equality.
\end{proof}

Here we introduce an example of the preceding lemma.
\begin{Example}
Similar to Example~\ref{groupnormal}, let~$G$ be a~group acting normally on a~set~$X$ from the right and
$H=\mathbb{C}G$ be the group algebra of~$G$ acting on $A=\operatorname{Fun}(X,{\mathbb C})$, the commutative algebra of all complex
valued functions on~$X$, by $(g f)(x)=f(x g)$.
It is easy to check that~$A$ is a~left~$H$-module algebra and furthermore that this action is commutative.
\end{Example}

One notes that for any bicrossed product Hopf algebra $H={\cal F} \acl {\cal U}$, the Hopf algebra ${\cal F}$ is
an~$H$-module algebra by the following action
\begin{gather}
\label{com5}
(f\acl u) g:= \varepsilon(f) (u g),
\qquad
f,g\in {\cal F}, \quad u\in {\cal U}.
\end{gather}

\begin{Example}
\label{cocom-action2}
Let $H:={\cal F} \acl {\cal U}$ be any bicrossed product Hopf algebra where ${\cal U}$ is a~cocommutative Hopf
algebra.
Then the following computation proves that the action in~\eqref{com5} is cocommutative
\begin{gather*}
(f\acl u)\ps{1}v_1\otimes (f\acl u)\ps{2}v_2= \big(f\ps{1}\acl u\ps{1}\ns{0}\big)v_1\otimes\big(f\ps{2}u\ps{1}\ns{1}\acl u\ps{2}\big)v_2
\\
 = \varepsilon\big(f\ps{1}\big)u\ps{1}\ns{0}v_1\otimes\varepsilon\big(f\ps{2}u\ps{1}\ns{1}\big)u\ps{2}v_2
 =\varepsilon(f)u\ps{1}v_1\otimes u\ps{2}v_2
 = (f\acl u)\ps{2}v_1\otimes (f\acl u)\ps{1}v_2.
\end{gather*}
Therefore any module~$M$ over~$H$, with the trivial coaction of~$H$, def\/ines a~$_{\cal F} H$-HCC.
\end{Example}

Here we introduce an example which shows that the categories of $_AH$-HCC and $_AH$-SAYD modules are dif\/ferent.
In fact we introduce an example of an $_AH$-HCC module which is not a~$_AH$-SAYD.
\begin{Example}
Consider ${\mathcal{H}}=\mathcal{H}_1^{\rm cop}\cong \mathcal{F}\acl \mathcal{U}$ and let $M={\mathcal{H}}$ be a~right
${\mathcal{H}}$-module via multiplication and left comodule by trivial coaction.
We have seen that the Hopf algebra ${\cal F}$ is a~left ${\mathcal{H}}$-module algebra by the action def\/ined
in~\eqref{com5} and this action is a~cocommutative action by Example~\ref{cocom-action2} and therefore that~$M$ is
a~$_{\cal F} {\mathcal{H}}$-HCC module.
We show that~$M$ is not a~$_{\cal F}{\mathcal{H}}$-AYD.
First one notices that using the formula
\begin{gather*}
S^{-1}(h)= \delta\big(S\big(h\ps{3}\big)\big)\delta\big(h\ps{1}\big)S\big(h\ps{2}\big),
\end{gather*}
we obtain
\begin{gather*}
 S^{-1}(\delta_1\acl 1)=- \delta_1\acl 1,
\qquad
  S^{-1}(\delta_2\acl 1)=- \delta_2\acl 1,
\\
  S^{-1}(\delta_1\acl Y)= \delta_1 \acl Y + \delta_1\acl 1,
\qquad
  S^{-1}(1\acl Y)= - 1 \acl Y,
\\
 S^{-1}(1\acl X)= -(1\acl X)+ (\delta_1\acl Y),
\\
  S^{-1}\big(\delta_1^2\acl Y\big)= -S^{-1}\big(2\delta_1^2\acl 1\big)-\big(\delta_1^2\acl Y\big),
\\
 S^{-1}(\delta_1\acl X)= (\delta_1\acl X)- (\delta_2\acl 1)-\big(\delta_1^2\acl Y\big).
\end{gather*}
Using~\eqref{fu} and the triviality of the coaction, by substituting $h= \delta_1\acl X$, $m= 1\acl 1$ and $a= \delta_1$
into the $_F{\mathcal{H}}$-AYD condition we have
\begin{gather*}
 S^{-1}((\delta_1\acl X)\ps{1}) (\delta_1\acl X)\ps{3} \delta_1\otimes m (\delta_1\acl X)\ps{2}
\\
\qquad
 =S^{-1}\big(\delta_1\ps{1}\acl X\ps{1}\ns{0}\big)\big(\delta_1\ps{3}X\ps{1}\ns{2}X\ps{2}\ns{1}\acl X\ps{3}\big)
 \delta_1\otimes (1\acl 1)\big(\delta_1\ps{2}X\ps{1}\ns{1}\acl X\ps{2}\ns{0}\big)
\\
\qquad
 =S^{-1}\big(\delta_1\ps{1}\acl X\ps{1}\ns{0}\big)X\ps{3} \delta_1\otimes \big(\delta_1\ps{2}X\ps{1}\ns{1}\acl X\ps{2}\big)
\\
\qquad
 =S^{-1}\big(\delta_1\ps{1}\acl X\ns{0}\big) \delta_1\otimes \big(\delta_1\ps{2}X\ns{1}\acl 1\big)+S^{-1}\big(\delta_1\ps{1}\acl 1\big)
\delta_1\otimes \big(\delta_1\ps{2}\acl X\big)
\\
\qquad\phantom{=}
 +S^{-1}\big(\delta_1\ps{1}\acl 1\big)X \delta_1\otimes \big(\delta_1\ps{2}\acl 1\big)
\\
\qquad
 =S^{-1}\big(\delta_1\ps{1}\acl X\big) \delta_1\otimes \big(\delta_1\ps{2}\acl 1\big)+S^{-1}\big(\delta_1\ps{1}\acl Y\big) \delta_1\otimes
\big(\delta_1\ps{2}\delta_1\acl 1\big)
\\
\qquad\phantom{=}
 +S^{-1}\big(\delta_1\ps{1}\acl 1\big) \delta_1\otimes \big(\delta_1\ps{2}\acl X\big)+ S^{-1}\big(\delta_1\ps{1}\acl 1\big)X \delta_1\otimes
\big(\delta_1\ps{2}\acl 1\big)
\\
\qquad
 =S^{-1}(\delta_1\acl X) \delta_1\otimes (1\acl 1) +S^{-1}(1\acl X) \delta_1\otimes (\delta_1\acl 1)
\\
\qquad\phantom{=}
 +S^{-1}(\delta_1\acl Y) \delta_1\otimes (1\delta_1\acl 1)+S^{-1}(1\acl Y) \delta_1\otimes (\delta_1\delta_1\acl 1)
\\
\qquad\phantom{=}
 +S^{-1}(\delta_1\acl 1) \delta_1\otimes (1\acl X)+ S^{-1}(1\acl 1) \delta_1\otimes (\delta_1\acl X)
\\
\qquad\phantom{=}
 +S^{-1}(\delta_1\acl 1)X \delta_1\otimes (1\acl 1)+S^{-1}(1\acl 1)X \delta_1\otimes (\delta_1\acl 1)
\\
\qquad
 =((\delta_1 \acl X)-(\delta_2 \acl 1)-(\delta_1^2\acl Y))\delta_1\otimes(1\acl 1)
 +(-(1\acl X)+(\delta_1\acl Y)) \delta_1
\\
\qquad\phantom{=}
 \otimes(\delta_1\acl 1)+((\delta_1\acl Y)+ \delta_1 \acl 1
) \delta_1\otimes (\delta_1\acl 1)-(1\acl Y)\delta_1\otimes (\delta_1^2\acl 1)
\\
\qquad\phantom{=}
 -(\delta_1\acl 1) \delta_1\otimes (1\acl X)+ (1\acl 1) \delta_1\otimes (\delta_1\acl X)
\\
\qquad\phantom{=}
 -(\delta_1\acl 1)X \delta_1\otimes (1\acl 1)+(1\acl 1)X \delta_1\otimes (\delta_1\acl 1)
\\
\qquad
 =-(X \delta_1)\otimes (\delta_1\acl 1)-(Y\delta_1)\otimes(\delta_1^2\acl 1)+\delta_1\otimes(\delta_1\acl X)+(X\delta_1)\otimes (\delta_1\acl 1)
\\
\qquad
 =\delta_1\otimes (\delta_1\acl X) -\delta_1 \otimes (\delta_1^2\acl 1)\\
\qquad
 \neq \delta_1\otimes (\delta_1\acl X)= \delta_1\otimes (1\acl 1)(\delta_1\acl X).
\end{gather*}
\end{Example}

\section{Cup products in Hopf cyclic cohomology}

In this section we show that all features of Hopf cyclic cohomology work well with the new coef\/f\/icients which we have
def\/ined in Section~\ref{S1}.
One of the most important features of Hopf cyclic cohomology is its cup product~\cite{atabey, kr5, rangipour} as
a~generalization of the Connes--Moscovici characteristic map~\cite{ConMos:HopfCyc}.
We show that this cup product works well with the new coef\/f\/icients in this paper.
One may construct an AYD module over a~Hopf algebra as the tensor product of a~YD module with a~AYD
module~\cite{HaKhRaSo1, staic}.
This shows that the category of AYD modules is a~${\cal C}$-category over the category of YD-modules.
In this section we show the same expectation for generalized modules is satisf\/ied.
In other words we prove that the category of $_CH$-SAYD modules is a~$\mathcal{C}$-category over the category of
$_CH$-YD modules.
\begin{Definition}
Let~$C$ be a~left~$H$-module coalgebra.
A~left-right module-comodule~$M$ over~$H$ is called an $_CH$-YD module if for $m\in M$, $h\in H$, and $c\in C$ we have
\begin{gather*}
(mh)\ns{-1} c \otimes (mh)\ns{0}= S^{-1}\big(h\ps{3}\big)m\ns{-1}h\ps{1}c \otimes m\ns{0}h\ps{2}.
\end{gather*}
\end{Definition}
\begin{Lemma}
Let~$M$ and~$N$ be right-left $_CH$-anti-Yetter--Drinfeld and Yetter--Drinfeld module respectively.
Then $M\otimes N$ is an $_CH$-anti-Yetter--Drinfeld module via
\begin{gather*}
(m\otimes n)h= mh\ps{2}\otimes nh\ps{1},
\qquad
m\otimes n\longmapsto m\ns{-1}n\ns{-1}\otimes m\ns{0}\otimes n\ns{-1}.
\end{gather*}
\end{Lemma}
\begin{proof}
The following computation shows that $M\otimes N$ is an $_CH$-anti-Yetter--Drinfeld module
\begin{gather*}
 ((m\otimes n)h )\ns{-1} c\otimes  ((m\otimes n)h )\ns{0}
 =\big(mh\ps{2}\otimes nh\ps{1}\big)\ns{-1} c\otimes \big(mh\ps{2}\otimes nh\ps{1}\big)\ns{0}
\\
\qquad
 =\big(mh\ps{2}\big)\ns{-1}(nh\ps{1})\ns{-1} c \otimes \big(mh\ps{2}\big)\ns{0}\otimes \big(nh\ps{1}\big)\ns{0}
\\
\qquad
 =S\big(h\ps{6}\big)m\ns{-1}h\ps{4}S^{-1}\big(h\ps{3}\big)n\ns{-1}h\ps{1} c\otimes m\ns{0}h\ps{5}\otimes n\ns{0}h\ps{2}
\\
\qquad
 =S\big(h\ps{4}\big)m\ns{-1}n\ns{-1}h\ps{1}\otimes m\ns{0}h\ps{3}\otimes n\ns{0}h\ps{2}
\\
\qquad
 =S\big(h\ps{3}\big)(m\otimes n)\ns{-1} h\ps{1}\otimes (m\otimes n)\ns{0} h\ps{2}.\tag*{\qed}
\end{gather*}
\renewcommand{\qed}{}
\end{proof}

The rest of this section is devoted to the cup product in Hopf cyclic cohomology with gene\-ra\-li\-zed coef\/f\/icients.

Let~$C$ be an~$H$-left module coalgebra acting from left on an~$H$-left module algebra~$A$ satisfying the following
conditions
\begin{gather}
\label{condition1}
  (hc)a=h(ca),
\qquad
 c(ab)=\big(c\ps{1}a\big)\big(c\ps{2} b\big),
\qquad
 c1= \varepsilon(c)1.
\end{gather}

Let $B= \Hom_{H}(C,A)$ denotes the convolution algebra of all~$H$-linear maps from~$C$ to~$A$ with the following
multiplication
\begin{gather*}
(f\ast g)(c)= f\big(c\ps{1}\big)g\big(c\ps{2}\big).
\end{gather*}
This algebra has the unit element $\eta \circ \varepsilon$ where $\eta: \mathbb{C}\longrightarrow A$ is the unit of~$A$.
Let $C^{*}_{H}(C, M)$ and $C^{*}_{H}(A, M)$ be the cocyclic modules which are def\/ined in~\eqref{module coalgebra}
and~\eqref{module algebra} where~$M$'s are right-left~$_CH$ and $_AH$-SAYD modules respectively.
Consider the diagonal complex
\begin{gather*}
C^{n,n}_{a,c}= C^{n}_{H}(C, M)\otimes C^{n}_{H}(A, M) =\Hom_{H}\big({M\otimes A^{\otimes n+1}}, \mathbb{C}\big)
\otimes
\big(M\otimes_{H}C^{\otimes (n+1)}\big),
\end{gather*}
which is a~cocyclic module by the following coface, codegeneracy and cyclic maps, $(d_n\otimes \delta_n, s_n\otimes
\sigma_n, t_n\otimes \tau_n)$.
For all $\varphi \in \Hom_{H}(M\otimes {A^{\otimes n+1}}, \mathbb{C})$, $f_i\in B$ and $m\otimes \widetilde{c}\in
M\otimes_{H}C^{\otimes (n+1)}$, we def\/ine the following map
\begin{gather*}
\nonumber
 \Psi_{a,c}:  \ C^{n,n}_{a,c}\longrightarrow \Hom\big(B^{\otimes (n+1)}, \mathbb{C}\big),
\\
\hphantom{\Psi_{a,c}:  \ }~\Psi_{a,c}(\varphi \otimes m\otimes c_0\otimes \dots \otimes c_n)(f_0\otimes \dots \otimes f_n)= \varphi(m \otimes_{H}
f_0(c_0)\otimes \dots \otimes f_n(c_n)).
\end{gather*}
Similar to~\cite{rangipour} one has the following statement.
\begin{Lemma}
The map $\Psi_{a,c}$ is a~well-defined map of cocyclic modules $C^{*,*}$ and $C^*(B)$.
\end{Lemma}

Using~\eqref{condition1}, 
one shows that the following map is an unital
algebra map
\begin{gather*}
 \chi: \ A\longrightarrow \Hom_{H}(C,A),
\\
 \hphantom{\chi: \  }~\chi(a)(c)=c(a).
\end{gather*}
Therefore we obtain a~cyclic map $\chi: C^*(B, \mathbb{C})\longrightarrow C^*(A, \mathbb{C})$.
As a~result we have the following cyclic map
\begin{gather*}
\Psi=\chi \circ \Psi_{a,c}: \ C^{n,n}_{a,c}\longrightarrow C^n(A, \mathbb{C}).
\end{gather*}

Now we def\/ine
\begin{gather*}
C^{p,q}_{a,c}= \bigoplus_{p+q=n} C_{H}^p(A, M) \otimes C_{H}^q(C, M),
\end{gather*}
to be a~tensor product of cocyclic modules which has a~mixed complex structure~\cite{AC0,Loday}.
The cyclic cohomology of this mixed complex is
\begin{gather*}
\bigoplus_{p,q} HC_{H}^p(A, M)\otimes HC_{H}^q(C, M).
\end{gather*}
As in~\cite{rangipour}, one uses the Alexander--Whitney map to obtain the following map
\begin{gather*}
{\rm AW}: \ HC_{H}^p(A, M)\otimes HC_{H}^q(C, M)\longrightarrow HC^{p+q}(C^{p,q}_{a,c})\longrightarrow HC^n(C^{n,n}_{a,c}).
\end{gather*}
Now we obtain the following statement.
\begin{Proposition}
Let~$C$ be an~$H$-left module coalgebra acting from the left on an~$H$-left mo\-du\-le algebra~$A$
satisfying~\eqref{condition1}. 
Let~$M$ be an $_AH $ and $_CH$-{\rm SAYD} module, and let $C^{*}_{H}(C, M)$ and $C^{*}_{H}(A, M)$ be the cocyclic modules
which are defined in~\eqref{module coalgebra} and~\eqref{module algebra}.
The following maps define cup products on the level of Hopf cyclic cohomology
\begin{gather*}
\sqcup=\Psi \circ {\rm AW}: \ HC_{H}^p(A, M)\otimes HC_{H}^q(C, M)\longrightarrow HC^{p+q}(A).
\end{gather*}
\end{Proposition}
\begin{proof}
This can be proved similarly as the author did in~\cite{rangipour}.
\end{proof}

\subsection*{Acknowledgements}
The authors of the manuscript are thankful to the organizers of Focus Program on Noncommutative Geometry and Quantum
Groups, which was held at Fields Institute June 3--28, 2013 for the invitation and the support.
Special thanks to P.M.~Hajac for his valuable comments and his unique attention to Hopf cyclic cohomology.
Last but not least, we would like to thank the referees for their extremely helpful comments.
This work is part of the project supported by the NCN grant 2011/01/B/ST1/06474.

\pdfbookmark[1]{References}{ref}
\LastPageEnding

\end{document}